\documentclass[smallextended,numbook,runningheads]{svjour3}
\usepackage{amsmath}
\smartqed  % flush right qed marks, e.g. at end of proof

\usepackage{amsfonts,amssymb}
\usepackage{epsfig}
\usepackage{booktabs}
\usepackage{xcolor}

\journalname{}
\date{ \phantom{b} \vspace{45mm}\phantom{e}}
\def\makeheadbox{\hspace{-4mm}Version of 24 April 2023}%, file: adiabatic.tex}

\def\R{\mathbb{R}}
 
\def\eps{\varepsilon}
\def\d{\mathrm{d}}
\def\e{{\mathrm e}}

\def\bigo{{\cal{O}}}
\def\wt{\widetilde}

\newcommand{\sfrac}[2]{\mbox{\footnotesize$\displaystyle\frac{#1}{#2}$}}
\def\half{\sfrac 12}

\def\bfe{\mathbf{e}}
\def\bfp{\mathbf{p}}
\def\bfu{\mathbf{u}}
\def\bfx{\mathbf{x}}
\def\bfA{\mathbf{A}}
\def\bfC{\mathbf{C}}
\def\bfF{\mathbf{F}}
\def\bfI{\mathbf{I}}
\def\bfK{\mathbf{K}}
\def\bfL{\mathbf{L}}
\def\bfM{\mathbf{M}}
\def\bfZ{\mathbf{Z}}
\def\bfzero{\mathbf{0}}

\newdimen\GGGlength
\newdimen\GGGheight
\newbox\GGGbox
\def\GGGput[#1,#2](#3,#4)#5{%
  \setbox\GGGbox\vbox{\hbox{#5}\kern0pt}%
  \GGGlength\wd\GGGbox%
  \divide\GGGlength by100 \multiply\GGGlength by#1%
  \GGGheight\ht\GGGbox%
  \divide\GGGheight by100 \multiply\GGGheight by#2%
  \put(#3,#4){\kern-\GGGlength\raise-\GGGheight\box\GGGbox}}

\begin{document}

\title{Leapfrog methods for relativistic charged-particle dynamics}

\titlerunning{Relativistic charged-particle dynamics}

\author{Ernst Hairer$^1$, Christian Lubich$^2$, Yanyan~Shi$^{2}$}
\authorrunning{E.\ Hairer, Ch.\ Lubich, Y.\ Shi}

\institute{
$^1$~Dept.\ de Math{\'e}matiques, Univ.\ de Gen{\`e}ve,
CH-1211 Gen{\`e}ve 24, Switzerland.\\
\phantom{$^1$~}\email{Ernst.Hairer@unige.ch}\\
$^2$~Mathematisches Institut, Univ.\ T\"ubingen, D-72076 T\"ubingen, Germany.\\
\phantom{$^2$~}\email{\{Lubich, Shi\}@na.uni-tuebingen.de}
%\\
%$^3$~LSEC, Academy of Mathematics and Systems Science, Chinese Academy of Sciences,\\ 
%\phantom{$^3$~}Beijing 100190, China; University of Chinese Academy of Sciences, Beijing 100049, China.\\
%\phantom{$^3$~}\email{shiyanyan1995@lsec.cc.ac.cn}
}

\date{ }

\maketitle

\begin{abstract}
A basic leapfrog integrator and its energy-preserving and variational / symplectic variants are proposed and studied for the numerical integration of the equations of motion of relativistic charged particles in an electromagnetic field. The methods are based on a four-dimensional formulation of the equations of motion. Structure-preserving properties of the numerical methods are analysed, in particular conservation and long-time near-conservation of energy and mass shell as well as preservation of volume in phase space. In the non-relativistic limit, the considered methods reduce to the Boris algorithm for non-relativistic charged-particle dynamics and its energy-preserving and variational / symplectic variants.
\bigskip

\bigskip

\end{abstract}

\section{Introduction}

In this paper we present and study numerical integrators for relativistic charged-particle dynamics. They are of interest in the context of particle methods for plasma physics when relativistic effects need to be taken into account.
The methods put forward here are based on a four-dimensional (4D) formulation of the equations of motion. The proposed leapfrog integrator and its energy-conserving and variational variants can be viewed as extensions of the standard Boris algorithm \cite{boris70rps} for non-relativistic charged-particle dynamics, which itself extends the leapfrog method (or St\"ormer--Verlet method) for particle motion in an electrostatic (potential-force) field; see e.g. Verlet~\cite{verlet67ceo} and \cite{hairer03gni}. However, the methods considered here differ substantially from the method proposed by Boris for the relativistic case, also given in \cite{boris70rps}. That method is based on a three-dimensional (3D) formulation of the equations of motion.

A concise review of the equations of motion of relativistic charged-particle dynamics in 3D and 4D formulations is given in Section 2.
In this paper we will work with the 4D formulation, which to us appears favourable for the numerical discretization. There is an ample literature on numerical integrators based on the 3D formulation, e.g. \cite{boris70rps,he16hov,higuera2017structure,ripperda2018comprehensive,vay2008simulation,zhang2015volume}, but only recently first numerical integrators based on the 4D formulation were proposed in the physical literature \cite{wang2016lorentz,matsuyama17hoi,xiao2019explicit,wang2021high}. Despite their algorithmic simplicity, the leapfrog methods of this paper have apparently not appeared in the literature on relativistic charged-particle dynamics before.

In Section 3 we propose and study the basic leapfrog integrator, which is explicit, or more precisely linearly implicit. It
preserves the mass shell and is also volume-preserving. Energy is in general not approximately conserved over long times along the numerical solutions, but often shows a random-walk behaviour, as it is the case also for the non-relativistic Boris algorithm \cite{hairer18ebo}. Energy {\it is} conserved in the special case of a constant electric field and is nearly conserved for all times in the case of a linear electric field, i.e., of a quadratic scalar potential. Unlike the Boris integrator in the non-relativistic setting \cite{hairer18ebo}, the method does not in general nearly preserve energy over long times in the case of a constant magnetic field, and not even for a vanishing magnetic field.

In Section 4 we consider a discrete-gradient implicit leapfrog integrator. It conserves energy and mass shell and is also volume-preserving.

In Section 5 we study a variational leapfrog integrator, which is also implicit. It exactly preserves a discrete energy that is close to the physical energy and it nearly preserves the mass shell over times that are exponentially long in the inverse step size. The method is symplectic and therefore conserves volume in the 8-dimensional phase space. The three different leapfrog methods of Sections 3, 4 and 5 coincide in the special case of constant electric and magnetic fields.

In Section 6 we discuss the non-relativistic limits of the three integrators. As it turns out, the relativistic explicit leapfrog method has the non-relativistic Boris method as its non-relativistic limit, the relativistic energy-preserving leapfrog method has a non-relativistic energy-preserving leapfrog method as its limit, and the relativistic variational leapfrog method has a known non-relativistic variational integrator as its limit. 

\section{Relativistic charged-particle dynamics}

In this preparatory section we  review the differential equations of motion together with their conserved quantities, the energy and the mass shell. They can be formulated in a 3D and a 4D setting, both of which appear in the physical literature;  see e.g. \cite[Section~8-4]{goldstein2011classical} and
\cite[Chapter~11]{jackson1999classical}.  Our presentation is partly inspired by \cite{wang2021high}. We will use the convention to denote 3D vectors and matrices by italic letters and 4D vectors and matrices by boldface letters.

% \subsection{Equations of motion and conserved quantities}\label{sec:model}
 \subsection{Three-dimensional formulation}
The equations of motion of a relativistic charged particle are usually formulated as
\begin{equation}\label{model}
\begin{aligned}
\frac{\d x}{\d t}&=\frac{u}{\gamma}\\
\frac{\d u}{\d t}&=\frac{u}{\gamma}\times B(x)+E(x)\\
 \gamma&=\sqrt{1+|u|^2}
\end{aligned}
\end{equation}
where $x(t)\in\mathbb R^3$ is the position at time $t$, $u(t)\in\mathbb R^3$ is the momentum, and $\gamma$ is the relativistic factor. Here, the physical units are chosen such that the speed of light is~1 
and the particle has mass~1 and electric charge~1. 
$E(x)$ and $B(x)$ are the given electric field and magnetic field, respectively, which can be expressed as $E(x)=-\nabla\phi(x)$ and $B(x)=\nabla\times A(x)$ with a scalar potential $\phi(x)\in\R$ and a vector potential $A(x)\in\R^3$.  We assume throughout the paper that $\phi$ and $A$ are arbitrarily differentiable (and hence also $E$ and $B$).

Equations \eqref{model} are the Euler--Lagrange equations $\frac d{dt}\nabla_v L = \nabla_x L$ for the Lagrangian
(with $v=u/\gamma$ and consequently $1/\gamma=\sqrt{1-|v|^2}$)
\[
L(x,v)=-\frac{1}{\gamma}+A(x)\cdot v-\phi(x).
\]
With the conjugate momenta 
$
p = \nabla_v L = \gamma v + A(x) = u +A(x),
$
the Hamiltonian $H=p\cdot v -L$ (energy of the charged particle) of \eqref{model} is 
(with $\gamma = \sqrt{1+|u|^2}=\sqrt{1+|p-A(x)|^2}$)
\begin{equation}\label{energy}
H(x,p)=\gamma+\phi(x).
\end{equation}
We will often write $H(x,\gamma)$ to emphasize the dependence of the energy on $x$ and $\gamma$ instead of the canonical variables  $(x,p)$.

\subsection{Four-dimensional formulation}
The numerical methods studied in this paper are based on a different formulation that works in 4D space-time.
By introducing the proper time $\tau$, the equations \eqref{model} can be expressed as 
\[
\begin{aligned}
%\left\{
\frac{\d t}{\d \tau}&=\gamma \qquad \frac{\d \gamma}{\d \tau}=E(x)\cdot u\\
\frac{\d x}{\d \tau}&=u \qquad
\frac{\d u}{\d \tau}=\gamma E(x)+u\times B(x),
\end{aligned}
\]
which is equivalent to the system of second-order differential equations
\begin{equation}\label{ode}
\begin{aligned}
\frac{\d^2 t}{\d \tau^2}&= E(x)\cdot \frac{\d {x}}{\d \tau}  \\
\frac{\d^2 x}{\d \tau^2}&=\, E(x)\ \frac{\d t}{\d \tau} - B(x) \times \frac{\d {x}}{\d \tau}\,.
\end{aligned}
\end{equation}
%With Lagrangian $L=\frac{1}{2}(-\dot{t}^2+|\dot{x}|^2)+A(x)\cdot\dot{x}-\dot{t}\phi(x)$, \eqref{ode} can be equivalently expressed as the Euler--Lagrangian system
%\[
%\frac{\d}{\d \tau}\left(\frac{\partial L}{\partial \dot{x}}\right)=\frac{\partial L}{\partial x}.
%\]
With the 4-dimensional position vector $\bfx=(t;x)=(t,x^\top)^\top$ and the matrix of the Minkowski metric $-dt^2+dx_1^2+dx_2^2+dx_3^2$,
$$
\bfM= \text{diag}(-1, 1, 1, 1),
$$
we can write \eqref{ode} in the compact form (with the dots standing for differentiation with respect to $\tau$),
\begin{equation}\label{ode2}
\bfM \ddot{\bfx}=\bfF(\bfx) \dot{\bfx} \quad\ \text{with skew-symmetric } \  \bfF=
\begin{pmatrix}
0&- E^\top\\
E&-\widehat{B}
\end{pmatrix},
\end{equation}
where $\widehat B$ is the $3\times 3$ matrix satisfying $\widehat{B} v = B \times v$ for all $v\in\R^3$. The matrix $\bfF$ represents the Faraday tensor. 

The first line of this 4-dimensional system of differential equations (with $\ddot t= \dot \gamma$) reads 
$-\dot \gamma = \nabla \phi(x)^\top \dot x \equiv \frac d{d\tau}\phi(x)$, which  imposes conservation of the energy \eqref{energy}.

The Lagrangian corresponding to \eqref{ode2} is (with $\bfu=(\gamma;u)$ as the variable representing $\dot \bfx$)
\begin{equation}\label{L4D}
\mathcal{L}(\bfx,\bfu)=\tfrac{1}{2}\,\bfu^\top  \bfM \bfu+\bfA(\bfx)^\top \bfu\quad\text{ with }\ \bfA(\bfx)=(-\phi(x);A(x)).
\end{equation}
%with the 4D vector potential 
%$$\bfA(\bfx)=(-\phi(x);A(x)).
%$$ 
Note that $\bfF=(\partial_\bfx \bfA)^\top  - \partial_\bfx \bfA$.

With  the conjugate momenta $\bfp=\nabla_{\bfu}\mathcal{L}=\bfM\bfu+\bfA(\bfx)$, the Hamiltonian (mass shell) is 
\begin{equation}\label{shell}
\mathcal{H}(\bfx,\bfp) =\tfrac12 \,\bfu^\top \! \bfM \bfu .
\end{equation}
We will often write $\mathcal{H}(\bfu)$ to emphasize the dependence of the mass shell on $\bfu$ instead of the canonical variables $(\bfx,\bfp)$.

Both Hamiltonians, the energy $H$ and the mass shell $\mathcal{H}$, are conserved along the solutions of the equivalent systems of differential equations \eqref{model} and \eqref{ode2}. As for every Hamiltonian system, the flow $\varphi_\tau: (\bfx(0),\bfp(0)) \mapsto (\bfx(\tau),\bfp(\tau))$ is symplectic and, in particular, preserves volume in phase space:
for every bounded open set $\Omega\subset \R^8$ and for every $\tau$ such that $\varphi_\tau(\bfx,\bfp)$ exists for all $(\bfx,\bfp)\in \Omega$, we have $\mathrm{vol} (\varphi_\tau(\Omega)) = \mathrm{vol} (\Omega)$. Since the transformation $(\bfx,\bfu)\mapsto(\bfx,\bfp)$ with $\bfp=\bfM\bfu+\bfA(\bfx)$ preserves volume, also the transformed flow $(\bfx(0),\bfu(0)) \mapsto (\bfx(\tau),\bfu(\tau))$ is volume-preserving.

Further conserved quantities arise if the 4D vector potential $\bfA$ is Lorentz-invariant: Suppose that for some matrix $\bfL\in \R^{4\times 4}$ with $\bfM \bfL + \bfL^\top \bfM=0$, the augmented vector potential $\bfA$ satisfies %then $\bigl(\e^{s\bfL}\bigr)^\top \bfM \e^{s\bfL} = \bfM$ for all real $s$ and
\begin{equation}\label{L-inv}
\bigl(\e^{s\bfL}\bigr)^\top \bfA(\e^{s\bfL}\bfx) =  \bfA(\bfx) \quad\text{ for all $s\in\R$ and all $\bfx\in \R^4$}.
\end{equation}
Then, the Lagrangian is invariant under the action of the group $(\e^{s\bfL})_{s\in\R}$, i.e. $\mathcal{L} (\e^{s\bfL}\bfx, \e^{s\bfL}\bfu)= \mathcal{L} (\bfx,\bfu)$. Therefore, Noether's theorem yields that
\begin{equation}\label{I-cq}
\mathcal{I}(\bfx,\bfp) = \bfp^\top \bfL\bfx, \qquad\text{where again }\ \bfp= \bfM\bfu+ \bfA(\bfx),
\end{equation}
is conserved along solutions of \eqref{ode2}, as is also verified by direct calculation. We will also write  $\mathcal{I}(\bfx,\bfu)$ to emphasize the dependence on $\bfx$ and $\bfu$ instead of the canonical variables $(\bfx,\bfp)$.

Finally we remark that the conservation of the mass shell \eqref{shell} and the invariant \eqref{I-cq} remain valid for a 4D vector
potential $\bfA(\bfx)=(-\phi(t,x);A(t,x))$ that also depends on time $t$. This is evidently not true for conservation of the energy~\eqref{energy}, but it remains true if only $A$ is time-dependent.

%\[
%\eta=\begin{pmatrix}
%-1&0&0&0\\
%0&1&0&0\\
%0&0&1&0\\
%0&0&0&1
%\end{pmatrix},
%\quad
%\]

%Since $\frac{\d t}{\d \tau}=\gamma=H_0-\phi(x)$, the above first equation can be written as
%\begin{equation}\label{secondorder}
%\frac{\d^2 x}{\d \tau^2}=\frac{\d {x}}{\d \tau}\times B(x)+ \widetilde{E}(x),
%\end{equation}
%where $\widetilde{E}(x)=(H_0-\phi(x))E(x)$.
%We assume $u_{\perp}=\bfp_{\perp}u=O(\eps)$, where $\bfp_{\perp}$ denotes the orthogonal projection onto the plane perpendicular to $B(x)$. 

%We are interested in the case of a strong magnetic field
%\[
%B(x)=B_\eps(\bfx)=\frac{1}{\eps}B_1(x), \quad 0\leq\eps\ll 1,
%\]
%where $B_1$ is smooth and independent of the small parameter $\eps$, with $|B_1(x)|\geq 1$ for all $x$.
%\section{Leapfrog methods based on the 4-dimensional formulation}
%%\begin{enumerate}
%%\item Leapfrog (LF)
%The numerical methods proposed and studied in this paper are based on the 4-dimensional differential equations \eqref{ode2}.
%We discuss two variants in this section: the basic explicit leapfrog method and an energy-preserving implicit leapfrog method.
%In a later section we will also consider a modified explicit leapfrog method for computing the guiding centre motion in the case of a strong non-uniform magnetic field using stepsizes that are much larger than the gyroperiod.

\section{The explicit leapfrog method based on the 4D formulation} 
\label{sec:lf}

\subsection{The explicit leapfrog method in two- and one-step form}
The simplest discretization of \eqref{ode2} replaces the derivatives by symmetric finite differences: Given a stepsize $h>0$ 
(with a fixed upper bound $h\le \bar h$) and approximations
$\bfx^n,\bfx^{n-1}$ to the augmented position vector $\bfx=(t;x)$ at $\tau=\tau_n=nh$ and $\tau=\tau_{n-1}$, the new position vector $\bfx^{n+1}$ is determined from
%The 4-dimensional position $\bfx=(t,x^\top)^\top$ is computed using the following two-step formulation 
\begin{equation}\label{LF}
\bfM\,\frac{\bfx^{n+1}-2\bfx^n+\bfx^{n-1}}{h^2}=\bfF(\bfx^n)\,\frac{\bfx^{n+1}-\bfx^{n-1}}{2h}.
\end{equation}
The augmented momentum $\bfu=(\gamma;u)=\dot \bfx$ is approximated by
\begin{equation}\label{Un}
\bfu^n=\frac{\bfx^{n+1}-\bfx^{n-1}}{2h}.
\end{equation}
%The 4-dimensional momentum $\bfu=(-\gamma;u)$ is approximated by
%\begin{equation}\label{\bfun}
%\bfu^n=\bfM\,\frac{\bfx^{n+1}-\bfx^{n-1}}{2h}.
%\end{equation}
For the actual computation a one-step formulation is preferable:
 Given the position and momentum approximation $(\bfx^n,\bfu^{n-1/2})$, the algorithm computes $(\bfx^{n+1},\bfu^{n+1/2})$ by leapfrog hopping:  
\begin{equation} \label{leapfrog}
\begin{aligned} 
\bfM\Bigl(\bfu^{n+\tfrac{1}{2}}-\bfu^{n-\tfrac{1}{2}}\Bigr)&=\frac{h}{2} \, \bfF(\bfx^n) \Bigl(\bfu^{n+\tfrac{1}{2}}+\bfu^{n-\tfrac{1}{2}}\Bigr)
\\
\bfx^{n+1}&=\bfx^n+h\bfu^{n+\tfrac{1}{2}}.
\end{aligned}
\end{equation}
Given the initial values $(\bfx^0,\bfu^0)$, the starting value $\bfu^{1/2}$ is chosen by first computing $u^{1/2}$ as the last three components of $\widetilde \bfu^{1/2}=\bfu^0 + \frac{h}{2}\bfM^{-1} \bfF(\bfx^0)\bfu^0$, then
$\gamma^{1/2}=\sqrt{1+| u^{1/2} |^2}$ and finally setting $\bfu^{1/2}=(\gamma^{1/2};u^{1/2})$.

In each step, the method requires evaluating the electric and magnetic fields $E$ and $B$ at the known current position $\bfx^n$ and then solving a 4-dimensional linear system with the matrix $\bfM-\tfrac12 h \bfF(\bfx^n)$.
The approximation \eqref{Un} of $\bfu$ at the grid points is obtained as $\bfu^n=\frac12\bigl(\bfu^{n+1/2}+\bfu^{n-1/2}\bigr)$.

While we claim no originality for this remarkably simple leapfrog method, it seems not to have been studied for relativistic charged-particle dynamics so far. This shall be done in this paper. We begin with a few observations. First, this is a symmetric (or time-reversible) method: interchanging the temporal superscripts $n+1$ and $n-1$ and replacing $h$ by $-h$ yields again the same method. As a consequence, it gives a second-order approximation as $h\to 0$: with $\tau^n=nh$,
$$
\| (\bfx^n,\bfu^n) - (\bfx( \tau^n),\bfu(\tau^n)) \| = \bigo(h^2) \qquad\text{for}\ nh \le \tau_{end},
$$
where the constant symbolized by $\bigo$ is independent of $h$ and $n$ with $nh\le \tau_{end}$, but depends on $\tau_{end}$.

\subsection{Preservation of mass shell and volume}

The following two theorems state important conservation properties.

\begin{theorem}\label{thm:mass-shell}
The leapfrog method \eqref{leapfrog} preserves the mass shell \eqref{shell}:
$$
\mathcal{H}\bigl(\bfu^{n+1/2}\bigr) = \mathcal{H}\bigl(\bfu^{n-1/2}\bigr).
$$
\end{theorem}

\begin{proof} Multiplying the first equation of \eqref{leapfrog} with $(\bfu^{n+1/2}+\bfu^{n-1/2})^\top$ and using the skew-symmetry of $\bfF(\bfx^n)$ yields
$$
\tfrac12 \bigl(\bfu^{n+1/2}\bigr)^\top \bfM \bfu^{n+1/2} - \tfrac12 \bigl(\bfu^{n-1/2}\bigr)^\top \bfM \bfu^{n-1/2} = 0,
$$
which is the stated result.
\qed
\end{proof}

As a consequence, for $\bfu^{n+1/2}=\bigl(\gamma^{n+1/2};u^{n+1/2}\bigr)$ the term $-\bigl(\gamma^{n+1/2})^2 +
\bigl|u^{n+1/2}\bigr|^2$ is independent of $n$, and hence
\begin{equation}\label{gamma-n}
\gamma^{n+1/2} = \sqrt{1+ \bigl|u^{n+1/2}\bigr|^2}\qquad\text{for all $n$,}
\end{equation}
provided that this relation holds true at $n=0$, as is ensured by our choice of starting value $\bfu^{1/2}=(\gamma^{1/2};u^{1/2})$.

\begin{theorem}\label{thm:volume}
The one-step map $(\bfx^n, \bfu^{n-1/2})\mapsto (\bfx^{n+1}, \bfu^{n+1/2})$ of the leapfrog method \eqref{leapfrog} is volume-preserving.
\end{theorem}

\begin{proof} 
The one-step map $\Phi_h:(\bfx,\bfu)\mapsto (\widehat \bfx,\widehat \bfu)$ given by
\begin{align*}
\widehat \bfx &= \bfx+ h \widehat \bfu \\
\widehat \bfu &= \mathrm{Cay}\bigl(\tfrac12 h \bfM^{-1} \bfF(\bfx)\bigr)\bfu
\end{align*}
with the Cayley transform $ \mathrm{Cay}(\bfZ)=(\bfI-\bfZ)^{-1}(\bfI+\bfZ)$ has the Jacobian
$$
D\Phi_h(\bfx,\bfu) = 
\begin{pmatrix} 
\bfI+ h \bfK \ & h\bfC \\ 
\bfK & \bfC
\end{pmatrix}
=
\begin{pmatrix} 
\bfI  \ & h\bfI \\ 
\bfzero & \bfI
\end{pmatrix}
\begin{pmatrix} 
\bfI & \bfzero \\ 
\bfK & \bfC
\end{pmatrix}
$$
with $\bfK=\partial \widehat \bfu / \partial \bfx$ and $\bfC= \partial \widehat \bfu / \partial \bfu = \mathrm{Cay}\bigl(\tfrac12 h \bfM^{-1} \bfF(\bfx)\bigr)$.
Since $\bfF(\bfx)$ is skew-symmetric and $\bfM$ is symmetric, the matrix $\bfZ= \tfrac12 h \bfM^{-1} \bfF(\bfx)$ trivially satisfies
$$
\bfM \bfZ+\bfZ^\top \bfM=0.
$$
This relation implies (see e.g.~\cite[Lemma IV.8.7]{hairer06gni}) that $\bfC=\mathrm{Cay}(\bfZ)$ satisfies
$$
\bfC^\top \! \bfM \bfC = \bfM
$$
and hence $|\det(\bfC)|=1$. This yields  $|\det \bigl( D\Phi_h(\bfx,\bfu) \bigr)|=1$ for all $(\bfx,\bfu)$, which shows that $\Phi_h$ is a volume-preserving map.
\qed
\end{proof}

\begin{remark} Since the map $(\bfx,\bfu)\mapsto (\bfx,\bfp)$ with $\bfp=\bfM\bfu+\bfA(\bfx)$ is volume-preserving, also the one-step map in the 8D phase space,
$(\bfx^n, \bfp^{n-1/2})\mapsto (\bfx^{n+1},\bfp^{n+1/2})$ with $\bfp^{n+1/2}=\bfM\bfu^{n+1/2} + \bfA(\bfx^{n+1})$, is volume-preserving.
\end{remark}

Long-time near-preservation of the energy $H$ of \eqref{energy} cannot be expected in general. However,
   in the special case of a quadratic electric potential $\phi$ the energy is nearly conserved over long times
   (see Theorem \ref{thm:mass-shell}), and for a linear potential it is exactly conserved (Remark~\ref{rem:energy}).
   This is different from the non-relativistic case, where the corresponding leapfrog method, the Boris method
   \cite{boris70rps}, does not preserve energy exactly in the case of a constant electric field.

\subsection{Energy behaviour: numerical experiments}

In the following examples we plot the relative error of the energy $H$ of \eqref{energy}
as a function of the proper time $\tau$. The time step size is $h$, the integration interval is
$[0,\tau_{end}]$, and the vertical axis covers the interval from $-ch^2$ to $+ch^2$, with $c$
depending on the particular example.

%\begin{example}[constant magnetic field]\label{example-2}
%\[
%\phi (x) = - \frac{0.01}{\sqrt{x_1^2 + x_2^2}}, \qquad B(x) = \bigl(0,0,1\bigr)^\top
%\]
%and initial values
%\begin{equation}\label{initial-2}
%x(0) = (0,1,0.1)^\top, \qquad u(0) = (0.09,0.05,0.2)^\top
%\end{equation}
%\end{example}

%\begin{figure}[t!]
%%\GGGinput{prog-re/}{boris-fig-2}
% \begin{picture}(0,0)
%  \epsfig{file=prog-re/boris-fig-2}
% \end{picture}%
%\begin{picture}(338.5,210.4)(  -0.7, 385.6)
%  \GGGput[  -30,  130](  0.60,594.66){$0.0015\, h^2$}
%  \GGGput[   50,  140](202.12,594.66){step size ~$h=0.1$}
%  \GGGput[  120,  -80](336.42,529.14){{\it Tend}$\,\, = 1\,000\,000$}
%  \GGGput[  -30,  130](  0.60,523.48){$0.0015\, h^2$}
%  \GGGput[   50,  140](202.12,523.48){step size ~$h=0.2$}
%  \GGGput[  120,  -80](336.42,458.07){{\it Tend}$\,\, = 1\,000\,000$}
%  \GGGput[  -30,  130](  0.60,452.41){$0.0015\, h^2$}
%  \GGGput[   50,  140](202.12,452.41){step size ~$h=0.4$}
%  \GGGput[  120,  -80](336.42,386.89){{\it Tend}$\,\, = 1\,000\,000$}
% \end{picture}
%\vspace{-3mm}
%\caption{Relative error of the Hamiltonian \eqref{energy} as a function of the
%proper time $\tau$ for the problem of Example~\ref{example-2}.
%\label{fig:boris-fig-2}}
%\end{figure}

\begin{figure}[t!]
%\GGGinput{prog-re/}{boris-fig-11}
 \begin{picture}(0,0)
  \epsfig{file=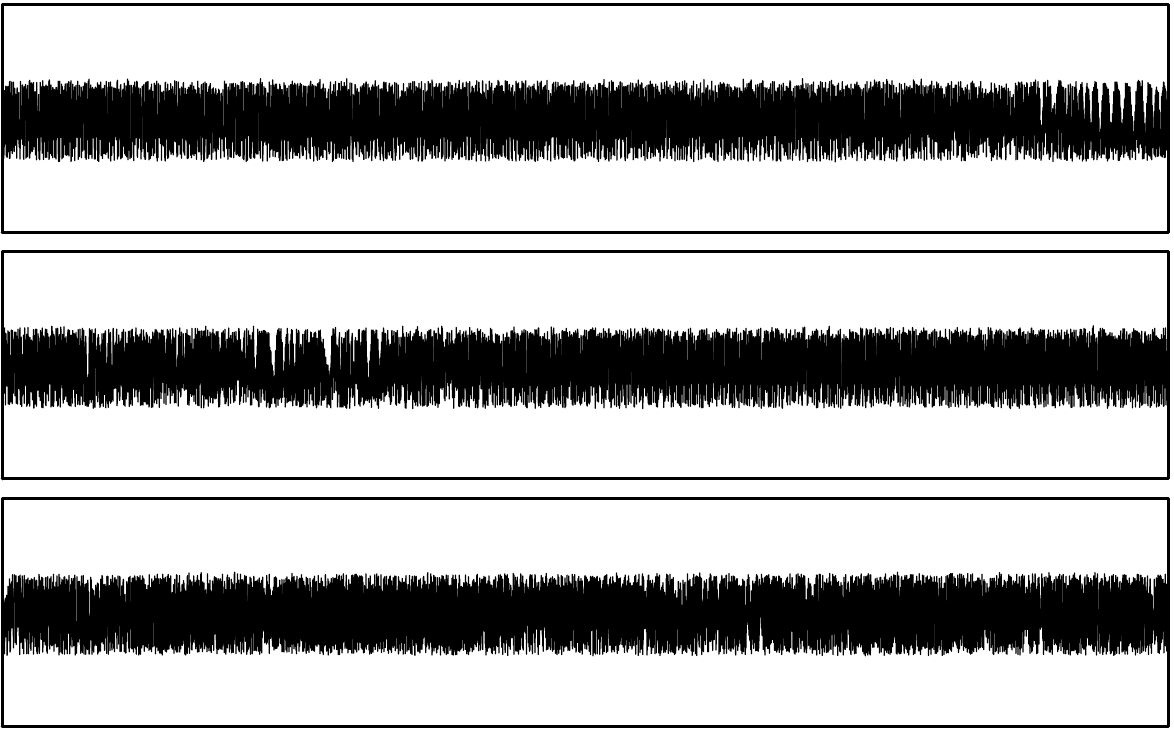}
 \end{picture}%
\begin{picture}(338.5,210.4)(  -0.7, 385.6)
  \GGGput[  -30,  130](  0.60,594.66){$2\, h^2$}
  \GGGput[   50,  140](202.12,594.66){step size ~$h=0.01$}
  \GGGput[  120,  -80](336.42,529.14){{\it Tend}$\,\, = 1\,000\,000$}
  \GGGput[  -30,  130](  0.60,523.48){$2\, h^2$}
  \GGGput[   50,  140](202.12,523.48){step size ~$h=0.02$}
  \GGGput[  120,  -80](336.42,458.07){{\it Tend}$\,\, = 1\,000\,000$}
  \GGGput[  -30,  130](  0.60,452.41){$2\, h^2$}
  \GGGput[   50,  140](202.12,452.41){step size ~$h=0.04$}
  \GGGput[  120,  -80](336.42,386.89){{\it Tend}$\,\, = 1\,000\,000$}
 \end{picture}
\vspace{-3mm}
\caption{Relative error of the energy \eqref{energy} as a function of the
proper time $\tau$ for the problem of Example~\ref{example-11}.
\label{fig:boris-fig-11}}
\end{figure}

\begin{example}[quadratic electric potential]\label{example-11}
\[
\phi (x) = x_1^2 +2x_2^2 +3x_3^2 - x_1, \qquad B(x) = \Bigl(0,0,\sqrt{x_1^2 + x_2^2}\,\Bigr)^\top
\]
with initial values $x(0) = (0,1,0.1)^\top$, $u(0) = (0.09,0.05,0.2)^\top$.
%\begin{equation}\label{initial-2}
%x(0) = (0,1,0.1)^\top, \qquad u(0) = (0.09,0.05,0.2)^\top.
%\end{equation}
\end{example}

%\begin{example}\label{example-20}
%\[
%\phi (x) = x_1^3 -x_2^3+\frac 15 x_1^4 + x_2^4 + x_3^4, \qquad B(x) = 
%\frac 12 \Bigl(x_2-x_3, x_1 + x_3, x_2-x_1\Bigr)^\top
%\]
%and initial values
%\begin{equation}\label{initial-20}
%x(0) = (0,1,0.1)^\top, \qquad u(0) = (0.09,0.55,0.3)^\top
%\end{equation}
%\end{example}
%
%
%\begin{figure}[t!]
%\centering
%%\GGGinput{prog-re/}{boris-fig-20}
% \begin{picture}(0,0)
%  \epsfig{file=prog-re/boris-fig-20}
% \end{picture}%
% \begin{picture}(338.5,210.4)(  -0.7, 385.6)
%  \GGGput[  -30,  130](  0.60,594.66){$8\,000\, h^2$}
%  \GGGput[   50,  140](202.12,594.66){step size ~$h=0.001$}
%  \GGGput[  120,  -80](336.42,529.14){{\it Tend}$\,\, = 100\,000$}
%  \GGGput[  -30,  130](  0.60,523.48){$8\,000\, h^2$}
%  \GGGput[   50,  140](202.12,523.48){step size ~$h=0.002$}
%  \GGGput[  120,  -80](336.42,458.07){{\it Tend}$\,\, = 100\,000$}
%  \GGGput[  -30,  130](  0.60,452.41){$8\,000\, h^2$}
%  \GGGput[   50,  140](202.12,452.41){step size ~$h=0.004$}
%  \GGGput[  120,  -80](336.42,386.89){{\it Tend}$\,\, = 100\,000$}
% \end{picture}
%\vspace{-3mm}
%\caption{Relative error of the Hamiltonian \eqref{energy} as a function of the
%proper time $\tau$ for the problem of Example~\ref{example-20}.
%\label{fig:boris-fig-20}}
%\end{figure}

\begin{figure}[t!]
\centering
%\GGGinput{prog-re/}{boris-fig-30}
 \begin{picture}(0,0)
  \epsfig{file=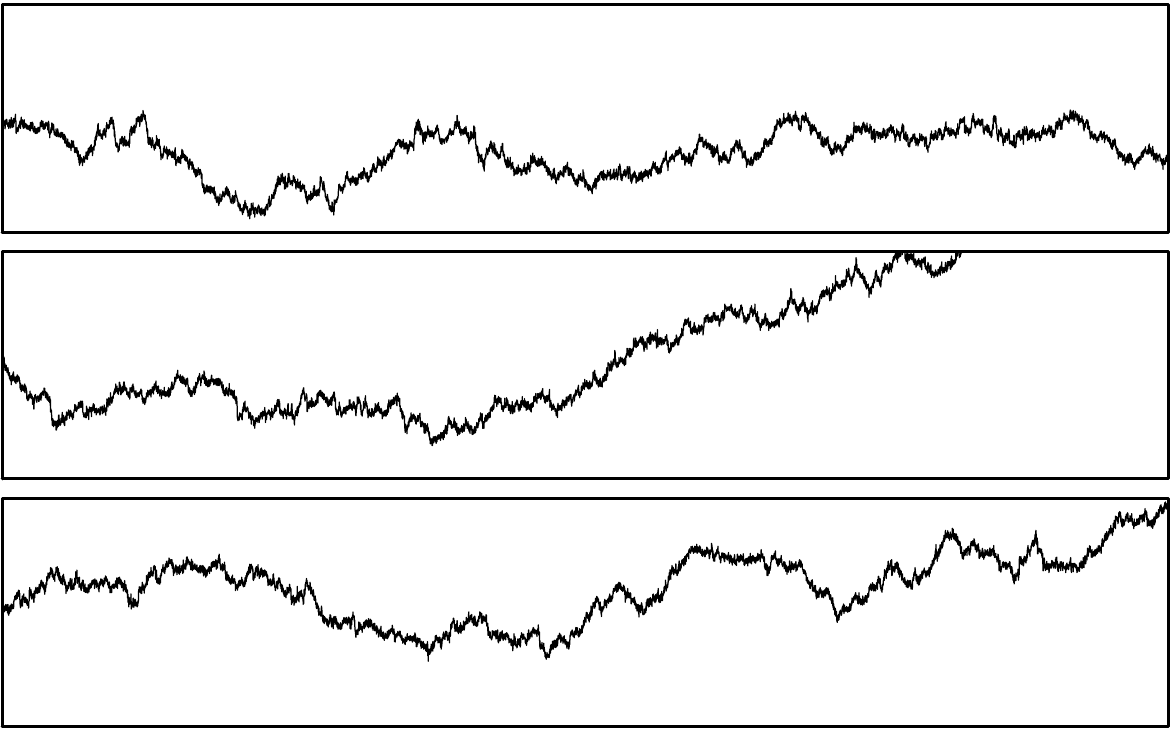}
 \end{picture}%
\begin{picture}(338.5,210.4)(  -0.7, 385.6)
  \GGGput[  -30,  130](  0.60,594.66){$4\,000\, h^2$}
  \GGGput[   50,  140](202.12,594.66){step size ~$h=0.0001$}
  \GGGput[  120,  -80](336.42,529.14){{\it Tend}$\,\, = 100\,000$}
  \GGGput[  -30,  130](  0.60,523.48){$4\,000\, h^2$}
  \GGGput[   50,  140](202.12,523.48){step size ~$h=0.0002$}
  \GGGput[  120,  -80](336.42,458.07){{\it Tend}$\,\, = 100\,000$}
  \GGGput[  -30,  130](  0.60,452.41){$4\,000\, h^2$}
  \GGGput[   50,  140](202.12,452.41){step size ~$h=0.0004$}
  \GGGput[  120,  -80](336.42,386.89){{\it Tend}$\,\, = 100\,000$}
 \end{picture}
\vspace{-3mm}
\caption{Relative error of the energy \eqref{energy} as a function of the
proper time $\tau$ for the problem of Example~\ref{example-30}.
\label{fig:boris-fig-30}}
\end{figure}

\begin{figure}[t!]
\centering
%\GGGinput{prog-re/}{boris-fig-21}
 \begin{picture}(0,0)
  \epsfig{file=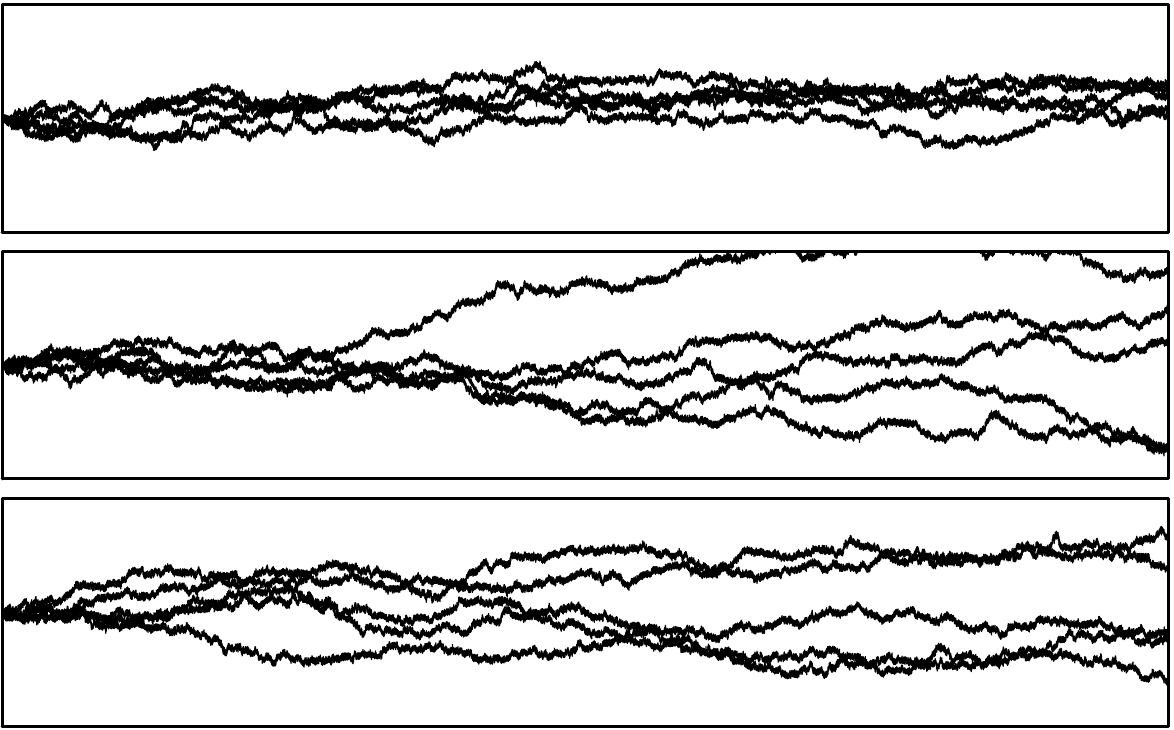}
 \end{picture}%
\begin{picture}(338.5,210.4)(  -0.7, 385.6)
  \GGGput[  -30,  130](  0.60,594.66){$5\,000\, h^2$}
  \GGGput[   50,  140](202.12,594.66){step size ~$h=0.001$}
  \GGGput[  120,  -80](336.42,529.14){{\it Tend}$\,\, = 100\,000$}
  \GGGput[  -30,  130](  0.60,523.48){$5\,000\, h^2$}
  \GGGput[   50,  140](202.12,523.48){step size ~$h=0.002$}
  \GGGput[  120,  -80](336.42,458.07){{\it Tend}$\,\, = 100\,000$}
  \GGGput[  -30,  130](  0.60,452.41){$5\,000\, h^2$}
  \GGGput[   50,  140](202.12,452.41){step size ~$h=0.004$}
  \GGGput[  120,  -80](336.42,386.89){{\it Tend}$\,\, = 100\,000$}
 \end{picture}
\vspace{-3mm}
\caption{Relative error of the energy \eqref{energy} as a function of the
proper time $\tau$ for the problem of Example~\ref{example-21}, for five initial values differing by less than $10^{-14}$.
\label{fig:boris-fig-21}}
\end{figure}

\begin{example}[non-quadratic electric potential]\label{example-30} 
\[
\phi (x) = x_1^3 -x_2^3+x_1^4/5 + x_2^4 + x_3^4, \qquad 
B(x) = \Bigl(0,0,\sqrt{x_1^2 + x_2^2}\,\Bigr)^\top
\]
with initial values $x(0) = (0,1,0.1)^\top$, $u(0) = (0.09,0.55,0.3)^\top$.
%\begin{equation}\label{initial-20}
%x(0) = (0,1,0.1)^\top, \qquad u(0) = (0.09,0.55,0.3)^\top
%\end{equation}
\end{example}

\begin{example}[constant magnetic field]\label{example-21}
\[
\phi (x) = x_1^3 -x_2^3+x_1^4/5 + x_2^4 + x_3^4, \qquad 
B(x) = \bigl(0,0,1\,\bigr)^\top
\]
and initial values as in Example~\ref{example-30}, but with perturbations $k\cdot 10^{-15}$ for $k=0,1,2,3,4$ added to the second component of $x(0)$.
\end{example}

In Figure~\ref{fig:boris-fig-11} we observe long-time near-conservation of energy in the case of a quadratic electric potential, but in Figure~\ref{fig:boris-fig-30} a random-walk behaviour is seen in the case of a non-quadratic electric potential. A random-walk behaviour is also observed for constant magnetic fields, see Figure~\ref{fig:boris-fig-21}, and even for the zero magnetic field (not shown).

The energy behaviour differs from that of the Boris method in the non-relativistic case: 
long-time near-conservation of energy was shown for the two cases of a constant magnetic field 
or a quadratic electric potential \cite[Theorem 2.1]{hairer18ebo}, and for the data of Example~\ref{example-30} a linear drift in the energy was observed in \cite[Example~5.1]{hairer18ebo}, whereas for a different (linear) magnetic field a random-walk behaviour was observed 
\cite[Example~5.2]{hairer18ebo}.

\subsection{Backward error  analysis}

The explicit leapfrog method can be written as
\begin{eqnarray*}
\frac{t^{n+1} - 2t^n + t^{n-1}}{h^2} & = & E(x^n)^\top \,\frac{x^{n+1} -x^{n-1}}{2h}\\[2mm]
\frac{x^{n+1} - 2x^n + x^{n-1}}{h^2} & = & E(x^n) \,\frac{t^{n+1} -t^{n-1}}{2h} -
\widehat B(x^n) \,\frac{x^{n+1} -x^{n-1}}{2h} ,
\end{eqnarray*}
and we let
\[ 
 \gamma^n = \frac{t^{n+1} -t^{n-1}}{2h} \qquad u^n = \frac{x^{n+1} -x^{n-1}}{2h}  .
\]
In the spirit of backward error analysis we search for an $h$-dependent modified differential equation
such that its solution $t(\tau ), x(\tau )$ formally satisfies $t(nh) = t^n$ and
$x(nh) = x^n$:
\begin{eqnarray*}
\frac{t(\tau +h) - 2t(\tau ) + t(\tau -h)}{h^2} & = & E\bigl( x(\tau )\bigr)^\top 
\,\frac{x(\tau +h) -x(\tau -h)}{2h}\\[2mm]
\frac{x(\tau +h) - 2x(\tau ) + x(\tau -h)}{h^2} & = & E\bigl( x(\tau )\bigr) \,\frac{t(\tau +h) -t(\tau -h)}{2h}
\\
&& - \widehat B \bigl( x(\tau )\bigr)  \,\frac{x(\tau +h) -x(\tau -h)}{2h} .
\end{eqnarray*}
A formal expansion into Taylor series yields (omitting the argument $\tau$)
\begin{eqnarray}
\ddot t + \frac{h^2}{12} \ddddot t+ \ldots  & = & E(x)^\top 
\,\Bigl( \dot x + \frac{h^2}6 \dddot x + \ldots \Bigr) \label{eq-t} \\[2mm]
\ddot x + \frac{h^2}{12} \ddddot x+ \ldots  & = & E(x) \,\Bigl( \dot t + \frac{h^2}6 \dddot t + \ldots \Bigr)
- \widehat B (x)  \,\Bigl( \dot x + \frac{h^2}6 \dddot x + \ldots \Bigr) . \qquad \label{eq-x}
\end{eqnarray}
With the $h$-dependent functions $\gamma (\tau ) = \bigl( t(\tau + h ) - t(\tau -h )\bigr) \big/ 2h$
and $u (\tau ) = \bigl( x(\tau + h ) - x(\tau -h )\bigr) \big/ 2h$
we also have
\[
\gamma = \dot t + \frac{h^2}6 \dddot t + \ldots ,\qquad 
u = \dot x + \frac{h^2}6 \dddot x + \ldots  .
\]
These relations can be formally inverted and give
\begin{equation}\label{tx-gu}
\dot t = \gamma - \frac{h^2}6 \ddot \gamma+ \ldots ,\qquad 
\dot x = u - \frac{h^2}6 \ddot u + \ldots  .
\end{equation}

\subsubsection{Quadratic electric potential}

With a quadratic scalar potential $\phi(x)$ we have long-time near-conservation of the energy
$H(x,\gamma)=\gamma+\phi(x)$ along numerical solutions obtained with the leapfrog method.

\begin{theorem}\label{thm:mass-shell}
Consider the leapfrog method \eqref{LF}--\eqref{Un} in the case of a quadratic potential $\phi (x) = \frac 12 \, x^\top Q x + q^\top x$. 
Provided that the numerical solution $(\bfx^n,\bfu^n)$ stays in a fixed compact set, 
the energy \eqref{energy} is conserved up to $\bigo(h^2)$,
$$
H(x^n,\gamma^n)=H(x^0,\gamma^0)+ \bigo(h^2),
$$
over long times $nh \le h^{-N}$ for arbitrary $N\ge 1$. The constant symbolized by $\bigo$ is independent of $n$ and $h$, but depends on the exponent $N$.
\end{theorem}

\begin{proof}
Inserting the first relation of \eqref{tx-gu} into \eqref{eq-t} yields
\begin{equation}\label{energymod}
\dot \gamma - \frac{h^2}{12} \dddot \gamma + \ldots =
E(x)^\top 
\,\Bigl( \dot x + \frac{h^2}6 \dddot x + \ldots \Bigr) .
\end{equation}
Along the solution of the modified differential equation
we have $E(x)^\top \dot x = - \frac{\d}{\d \tau} \phi (x)$, and we know from
\cite{hairer18ebo} that for a quadratic potential
$\phi (x)$,
the function $E(x)^\top x^{(k)}$  can be written as a total differential for odd values of $k$.
Consequently, there exist $h$-independent functions
$H_{2j}(x,\gamma )$ such that the function
\[
H_h(x,\gamma ) = H(x,\gamma ) + h^2 H_2(x,\gamma ) + h^4 H_4(x,\gamma ) 
+ \ldots ~ ,
\]
truncated at the $\bigo (h^N)$ term, satisfies
\[
\frac {\d} {\d \tau} H_h(x,\gamma ) = \bigo (h^N)
\]
along solutions of the modified differential equation. This yields the result as in \cite[Theorem 2.1]{hairer18ebo}.
\qed
\end{proof}

\subsubsection{Constant magnetic field}

In the non-relativistic case, the Boris method has long-time near-conservation of energy also for a constant magnetic field \cite[Theorem 2.1]{hairer18ebo}. The numerical experiments given above indicate that this is no longer true in the relativistic case with the leapfrog method~\eqref{LF}. Let us explain where the different behaviour comes from.
We write \eqref{energymod} as 
\begin{equation}\label{energyB}
\frac {\d}{\d \tau} \Bigl( \gamma  + \phi (x) \Bigr) - \frac{h^2}{12} \dddot \gamma + \ldots =
E(x)^\top 
\,\Bigl( \frac{h^2}6 \dddot x + \ldots \Bigr) .
\end{equation}
Computing $E(x)$ from \eqref{eq-x} and inserting the resulting expression into
\eqref{energyB} yields
\[
\Bigl( \dot t + \frac{h^2}6 \dddot t + \ldots \Bigr)^{-1}\biggl((\ddot x + \frac{h^2}{12} \ddddot x+ \ldots 
+ \widehat B (x)  \,\Bigl( \dot x + \frac{h^2}6 \dddot x + \ldots \Bigr) \biggr)^\top
\,\Bigl( \frac{h^2}6 \dddot x + \ldots \Bigr)
\]
for the right-hand side of \eqref{energyB}. If the magnetic field is constant,
it follows from \cite{hairer18ebo} that ${x^{(k)}}^\top \widehat B \, x^{(l)}$ is a total differential
when $k-l$ is an even integer. Consequently, there exist functions $F_{2j}(x,\dot x,\gamma )$
such that
\[
\frac {\d}{\d \tau} \Bigl( \gamma  + \phi (x) \Bigr) - \frac{h^2}{12} \dddot \gamma + \ldots =
\Bigl( \dot t + \frac{h^2}6 \dddot t + \ldots \Bigr)^{-1}\frac{\d}{\d \tau}
\Bigl(h^2 F_2(x,\dot x,\gamma )+ \ldots \Bigr) .
\]
Unfortunately, this relation does not permit us to prove long-time near-conser\-vation of energy for the explicit leapfrog method. 
%This is confirmed by the following example. This behaviour is different from non-relativistic charged-particle dynamics, where the corresponding leapfrog method, the Boris method, is known to preserve energy up to $\bigo(h^2)$ for very long times in the case of a constant magnetic field~\cite{hairer18ebo}.

\section{An energy-preserving implicit leapfrog method}

The explicit leapfrog method in the preceding subsection does not preserve the energy $H$ of \eqref{energy}. We describe an implicit energy-preserving variant that uses discrete gradients. Given the potential $\phi$, a continuous map $\overline\nabla\phi:\R^3\times\R^3\to\R^3$ is a discrete gradient of $\phi$ if the following two conditions are satisfied for all $x,\widehat x\in\R^3$:
\begin{align*}
\overline{\nabla}\phi(\widehat x,x) \cdot (\widehat x-x) = \phi(\widehat x)-\phi(x), \qquad  \overline{\nabla}\phi(x,x) = \nabla \phi(x).
\end{align*}
Well-known examples are the midpoint discrete gradient \cite{gonzalez96tia}
$$%\begin{equation}\label{dgg}
\overline \nabla \phi(\widehat x, x) \, = \,
\nabla \phi(\overline x) +
{ \phi(\widehat x) - \phi(x) - \nabla \phi(\overline x)^T \Delta x
\over \| \Delta x \|^2} \, \Delta x
$$%\end{equation}
with $\overline x = \tfrac12(\widehat x + x)$ and $\Delta x= \widehat x -x$ and the average vector field \cite{celledoni12per}
$$
\overline \nabla \phi(\widehat x, x) = \int_0^1 \nabla \phi(x+\theta(\widehat x -x))\,d\theta.
$$
To obtain an energy-conserving method, we replace \eqref{LF} by
\begin{equation}\label{ECLF}
\bfM\,\frac{\bfx^{n+1}-2\bfx^n+\bfx^{n-1}}{h^2}=\overline \bfF^n\,\frac{\bfx^{n+1}-\bfx^{n-1}}{2h},
\end{equation}
where we choose
$$\overline \bfF^n= \begin{pmatrix}
0&(\overline{\nabla}\phi^n)^\top\\
-\overline{\nabla}\phi^n &-\widehat{B}(x^n)
\end{pmatrix}
$$
with the discrete gradient %$\overline{\nabla}\phi^n=\overline{\nabla}\phi(x^{n+1/2},x^{n-1/2})$ for
$$%\begin{equation}\label{mid}
\overline{\nabla}\phi^n=\overline{\nabla}\phi(x^{n+1/2},x^{n-1/2}) \quad\text{for } \ x^{n\pm1/2}=\tfrac12(x^n+x^{n\pm1}).
$$%\end{equation}
The method has a one-step formulation analogous to \eqref{leapfrog}, with $\overline \bfF^n$ instead of $\bfF(\bfx^n)$. This yields augmented momentum approximations $\bfu^{n+1/2}$.
%The formulation analogous to \eqref{leapfrog} then becomes, with $\bfu^{n+1/2}=(\gamma^{n+1/2},u^{n+1/2})$,
%\begin{equation} \label{ec-leapfrog}
%\begin{aligned} 
%\bfM\Bigl(\bfu^{n+\tfrac{1}{2}}-\bfu^{n-\tfrac{1}{2}}\Bigr)&=\frac{h}{2} \, \overline \bfF^n \Bigl(\bfu^{n+\tfrac{1}{2}}+\bfu^{n-\tfrac{1}{2}}\Bigr)
%\\
%\bfx^{n+1}&=\bfx^n+h\bfu^{n+\tfrac{1}{2}}.
%\end{aligned}
%\end{equation}
In contrast to method \eqref{leapfrog}, this method is now implicit in $x^{n+1}$, since $\overline{\nabla}\phi^n$ depends on $x^{n+1}$. Also note that now $x^{n-1}$ enters $ \overline \bfF^n $, and so we do not have a pure one-step scheme as in the explicit leapfrog method \eqref{leapfrog}.

\begin{theorem}\label{thm:ec}
The discrete-gradient leapfrog method %\eqref{ec-leapfrog} 
preserves both the mass shell \eqref{shell} and the energy \eqref{energy}: for all $n$,
\begin{align*}
\mathcal{H}\bigl(\bfu^{n+1/2}\bigr) &= \mathcal{H}\bigl(\bfu^{1/2}\bigr)
\\
H(x^{n+1/2}, \gamma^{n+1/2}) &= H(x^{1/2}, \gamma^{1/2}).
\end{align*}
\end{theorem}

\begin{proof}
We note that
\[
\frac{u^{n+1/2}+u^{n-1/2}}{2}=u^n=\frac{x^{n+1}-x^{n-1}}{2h}=\frac{x^{n+1/2}-x^{n-1/2}}{h}.
\]
From the first line in the first equation of the one-step formulation \eqref{leapfrog} with $ \overline \bfF^n $ instead of $\bfF(x^n)$ and the definition of a discrete gradient, we therefore obtain
\[
\begin{aligned}
-\frac{\gamma^{n+1/2}-\gamma^{n-1/2}}{h}&=\overline{\nabla}\phi^n\cdot\frac{u^{n+1/2}+u^{n-1/2}}{2}\\
&=\overline{\nabla}\phi(x^{n+1/2},x^{n-1/2})\cdot\frac{x^{n+1/2}-x^{n-1/2}}{h}\\
&=\frac{\phi(x^{n+1/2})-\phi(x^{n-1/2})}{h}.
\end{aligned}
\]
In view of \eqref{energy}, this proves the energy conservation $H^{n+1/2}=H^{n-1/2}$.\\
Conservation of mass shell is proved in the same way as in Theorem~\ref{thm:mass-shell}. 
\qed
\end{proof}

\begin{remark} 
\label{rem:energy}
In the special case of a constant electric field $E$ the two integrators of this and the preceding subsection coincide, since then
$\phi(x)=-E^\top x$ and hence $-\overline{\nabla}\phi^n = - \nabla \phi(x^n) =E$. In particular, in this case the explicit leapfrog integrator of Section~\ref{sec:lf} conserves the energy $H$.
\end{remark}

The discrete-gradient leapfrog method cannot be written as a one-step method, because $\overline{\nabla}\phi^n$ depends on $x^{n+1},x^n,x^{n-1}$. However, we have the following result.

\begin{theorem}\label{thm:volume-dg}
With every $\bfx^n$ in an open subset of $\R^4$, let $\bfx^{n-1}$ be associated in a differentiable but otherwise arbitrary way.
Then, the one-step map $(\bfx^n, \bfu^{n-1/2})\mapsto (\bfx^{n+1}, \bfu^{n+1/2})$ of the discrete-gradient leapfrog method is volume-preserving for any choice of the discrete gradient.
\end{theorem}

\begin{proof} The proof follows the lines of the proof of Theorem~\ref{thm:volume}, noting that $\overline \bfF^n$ is still skew-symmetric and that for any matrix $\bfK=\partial \widehat \bfu / \partial \bfx$ the one-step map is volume-preserving.
\qed
\end{proof}

%\subsection{Numerical experiments}

\section{Variational leapfrog integrator}
The variational integrator described here is constructed in the same way as is done in the interpretation of the St\"ormer--Verlet--leapfrog method  as a variational integrator;
see e.g. \cite[Section 1.6]{hairer03gni}.  The integral of the Lagrangian \eqref{L4D} over a time step is approximated in two steps: the path $x(t)$ of positions is approximated by the linear interpolant of the endpoint positions, and the integral is approximated by the trapezoidal rule. This yields the following approximation to the action integral over a time step, the `discrete Lagrangian'
\begin{equation}\label{Lh}
\mathcal{L}_h(\bfx^n,\bfx^{n+1}) = \frac h2 \mathcal{L}(\bfx^n, \bfu^{n+1/2}) + \frac h2 \mathcal{L}(\bfx^{n+1}, \bfu^{n+1/2}) 
\end{equation}
with $\bfu^{n+1/2} =\bigl(\bfx^{n+1}-\bfx^{n}\bigr)/h$. Extremizing this expression leads to the discrete Euler-Lagrange equations, which are the equations of the following implicit two-step method:
\begin{align}\label{varint}
 &\bfM\,\frac{\bfx^{n+1}-2 \bfx^n + \bfx^{n-1}}{h^2} =
 \\
 \nonumber
 &\qquad {\bfA}'(\bfx^n)^\top \,\frac {\bfx^{n+1}-\bfx^{n-1}}{2h}
-  \frac{{\bfA}(\bfx^{n+1})- {\bfA}(\bfx^{n-1})}{2h},
\end{align}
with the derivative matrix ${\bfA}'(\bfx)=\partial_\bfx {\bfA}(\bfx)=\bigl( \partial_j {\bfA}_i(\bfx)\bigr)_{i,j=1}^4$ for $\bfx=(t;x)$; cf.~\cite{hairer20lta,hairer22lsi,webb14sio} for the analogous method in the non-relativistic case.

This method is again complemented with the augmented-momentum approximation~\eqref{Un}.  
The method can equivalently be written as a perturbation of the explicit leapfrog method:
\begin{align} \nonumber
 &\bfM\,\frac{\bfx^{n+1}-2 \bfx^n + \bfx^{n-1}}{h^2} = \bfF(\bfx^n) \,\frac {\bfx^{n+1}-\bfx^{n-1}}{2h}
 \\
\label{varint-leapfrog-like}
 &\qquad\qquad  + {\bfA}'(\bfx^n) \,\frac {\bfx^{n+1}-\bfx^{n-1}}{2h}
-  \frac{{\bfA}(\bfx^{n+1})-{\bfA}(\bfx^{n-1})}{2h}\, .
\end{align}
The method has a one-step formulation similar to \eqref{leapfrog} of the explicit leapfrog algorithm, adding the correction term of \eqref{varint-leapfrog-like} in the first line of \eqref{leapfrog}.
It is, however, an implicit method, because the vector potential $\bfA$ is evaluated at the new position $\bfx^{n+1}$.

\begin{remark}
The correction to the explicit leapfrog method as given in the second line of \eqref{varint-leapfrog-like} vanishes for linear $\bfA(\bfx)$. Hence, the variational integrator and the explicit leapfrog method as well as the energy-preserving leapfrog method coincide for constant fields $B$ and $E$.
\end{remark}

A different one-step formulation of \eqref{varint} is given by the map $(\bfx^n,\bfp^n)\mapsto(\bfx^{n+1},\bfp^{n+1})$ that is implicitly defined by
\begin{align}
\bfp^n &= - \partial_1 \mathcal{L}_h(\bfx^n,\bfx^{n+1}) \nonumber \\
&= \Bigl( \bfM - \tfrac12\, h \bfA'(\bfx^n)^\top \Bigr) \bfu^{n+1/2} + \half\Bigl( \bfA(\bfx^n)+\bfA(\bfx^{n+1}) \Bigr)
\label{P-n}
\\
\bfp^{n+1} &= \partial_2\mathcal{L}_h(\bfx^n,\bfx^{n+1})  \nonumber
\\
&= \Bigl( \bfM + \tfrac12\, h \bfA'(\bfx^{n+1})^\top \Bigr) \bfu^{n+1/2} + \half\Bigl( \bfA(\bfx^n)+\bfA(\bfx^{n+1}) \Bigr),
\nonumber%\label{P-n+1}
\end{align}
where again $\bfu^{n+1/2} =(\gamma^{n+1/2};u^{n+1/2}) = \bigl(\bfx^{n+1}-\bfx^{n}\bigr)/h$.
This map yields the numerical positions $\bfx^n$ of the discrete Euler--Lagrange equations \eqref{varint} and is known to be symplectic (and hence also volume-preserving); see Theorem~VI.6.1 in~\cite{hairer06gni}, which can be traced back to Maeda~\cite{maeda82lfo}.

In the continuous problem, the first component of  $\bfp=\bfM\bfu + \bfA(\bfx)$ equals $-\gamma-\phi(x)= -H(x,\gamma)$, i.e. the negative energy \eqref{energy}. 
Since the 4D vector potential $\bfA(\bfx)=(-\phi(x);A(x))$ is independent of $t$ when $\phi$ and $A$ do not depend on time, the first column of $\bfA'(\bfx)$ is zero. This implies that in the discrete equations \eqref{P-n}, the first component of both $\bfp^n$ and $\bfp^{n+1}$ equals the negative of
\begin{equation}\label{energy-h}
H^{n+1/2} := \gamma^{n+1/2}+\half \Bigl( \phi(x^n)+\phi(x^{n+1})\Bigr).
\end{equation}
So we have proved the following conservation result.

\begin{theorem} \label{thm:H-varint}
The variational leapfrog integrator preserves the discrete energy \eqref{energy-h}:
$$
H^{n+1/2}=H^{1/2}
\qquad\text{for all}\ n.
$$
\end{theorem}

As long as $|x^{n+1/2}|\le M_0$ and $|u^{n+1/2}|\le M_1$, Taylor expansion of $\phi(x)$ at $x^{n+1/2}$ using $x^{n+1/2 \pm 1/2} = x^{n+1/2} \pm \tfrac12 h u^{n+1/2}$ shows that
$$
H^{n+1/2} = H(x^{n+1/2},\gamma^{n+1/2}) + \bigo(h^2),
$$
and hence Theorem~\ref{thm:H-varint} yields near-conservation of the energy \eqref{energy} up to $\bigo(h^2)$:
\begin{equation}\label{varint-H-near}
H(x^{n+1/2},\gamma^{n+1/2}) = H(x^{1/2},\gamma^{1/2}) + \bigo(h^2),
\end{equation}
where the constant symbolized by $\bigo$ is independent of $h$ and $n$, but depends on the bounds $M_0$ and $M_1$.
 
For the mass-shell we have the following near-conservation result over exponentially long times.

\begin{theorem} \label{thm:Hcal-varint}
Assume that the potentials $\phi$ and $A$ are analytic in a domain $D\subset \R^4$.
Provided that the numerical solution $(\bfx^n,\bfu^{n+1/2})$ stays in a compact set, and $\bfx^n\in D$, 
the variational leapfrog integrator preserves the mass shell~\eqref{shell} up to $\bigo(h^2)$,
$$
\mathcal{H}(\bfu^{n+1/2}) = \mathcal{H}(\bfu^{1/2}) + \bigo(h^2),
$$
where the constant symbolized by $\bigo$ is independent of $n$ and $h$ 
over exponentially long times $nh \le e^{c/h}$ with $c>0$.
\end{theorem}

\begin{proof}
The result is obtained from the known theory of variational / symplectic integrators; see e.g. \cite[Chapters VI and IX]{hairer06gni}. Backward error analysis shows that symplectic methods, such as \eqref{P-n} in our case,
nearly conserve the Hamiltonian, which here is the mass shell \eqref{shell}, over very long times; see Theorem IX.8.1 in \cite{hairer06gni}, which goes back to Benettin \& Giorgilli~\cite{benettin94oth}. Under the given assumptions, this yields that
the mass shell~\eqref{shell} is nearly conserved up to $\bigo(h^2)$,
\begin{equation}\label{HXP}
\mathcal{H}(\bfx^n,\bfp^n) = \mathcal{H}(\bfx^0,\bfp^0) + \bigo(h^2)\quad\text{ for $nh \le e^{c/h}$.}
\end{equation}
Using \eqref{P-n},  a calculation shows that the mass shell expressed in terms of $\bfu$ satisfies along the actually computed numerical solution $(\bfx^n,\bfu^{n+1/2})$
$$
\mathcal{H}(\bfu^{n+1/2}) = \half \Bigl( \mathcal{H}(\bfx^n,\bfp^n) + \mathcal{H}(\bfx^{n+1},\bfp^{n+1}) \Bigr) + \bigo(h^2).
$$
The relation \eqref{HXP} therefore implies that
$
\mathcal{H}(\bfu^{n+1/2}) = \mathcal{H}(\bfu^{1/2}) + \bigo(h^2)
$
over exponentially long times $nh \le e^{c/h}$.
\qed
\end{proof}

\begin{remark}
The invariant $\mathcal{I}$ of \eqref{I-cq} in the case of a Lorentz-invariant 4D vector potential $\bfA(\bfx)$ satisfying \eqref{L-inv} is exactly preserved by the variational integrator:
$$
\mathcal{I}(\bfx^n,\bfp^n)=\mathcal{I}(\bfx^0,\bfp^0) \qquad\text{for all $n$}.
$$
This follows directly from the discrete Noether's theorem \cite[Theorem~VI.6.7]{hairer06gni}, which states that invariants resulting from Noether's theorem are preserved by variational / symplectic integrators that preserve the symmetry, as is the case for the linear group action $(\e^{s\bfL})_{s\in\R}$ with \eqref{L-inv}. 

The conservation of the first component of $\bfp$ along the numerical solution, and hence the energy conservation of Theorem~\ref{thm:H-varint}, can also be seen as an instance of the discrete Noether's theorem, in view of the invariance of the discrete Lagrangian under the group action $\bfx \mapsto \bfx+s\bfe_1$, $s\in\R$, with $\bfe_1$ the first 4-dimensional unit vector, which yields the preservation of $ \bfe_1^\top \bfp^n = -H^{n+1/2}$.
\end{remark}

\section{Non-relativistic limit}
In the case of a small momentum $u=\eps\wt u$ with $\eps\ll1$ and moderately bounded $\wt u$  we have
$$
\gamma = \sqrt{1+\eps^2|\wt u|^2} = 1 + \tfrac12 \eps^2 |\wt u|^2 + \bigo(\eps^4).
$$
With a small potential that scales as $\phi(x)=\eps^2 \wt\phi(x)$ with moderately bounded $\wt\phi$,  the energy is
$$
H(x,\gamma)=\gamma+\phi(x) = 1 + \eps^2\bigl( \tfrac12 |\wt u|^2 + \wt \phi \bigr) + \bigo(\eps^4),
$$
where the constant term 1 is irrelevant, and the second term on the right-hand side is the $\eps^2$-scaled non-relativistic energy. In the following we omit the tilde on the variables.

With the scaling $u \to \eps u$, $\phi\to \eps^2 \phi$, $A \to \eps A$, and the rescaling of time as $\eps t \to t$ and $\eps\tau \to \tau$, the equations of motion \eqref{ode2} become $\dot t = \gamma$, $\dot x = u$ and
$$
\begin{pmatrix}
- \dot \gamma \\ \dot u
\end{pmatrix}
= 
\begin{pmatrix}
- \eps^2 E(x)^\top u \\ E(x)\gamma - \widehat B(x) u
\end{pmatrix}.
$$
On an interval of length $\tau=o(\eps^{-2})$, and under the assumption of bounded $x(\tau)$ and (rescaled) $u(\tau)$ on this interval, we obtain
$$
\dot t = \gamma = 1 + \bigo(\eps^2 \tau).
$$
In the limit $\eps\to 0$, we thus obtain $\gamma=1$ and $t=\tau$ and the familiar classical equations of motion of a charged particle,
$$
\dot x = u, \qquad \dot u = E(x) + u \times B(x).
$$
We next consider what the corresponding limit equations are for the three leapfrog methods studied in this paper.
With the above scaling and the appropriate rescaling $\eps h \to h$ of the stepsize, the limit equation for the explicit leapfrog method
\eqref{LF} becomes the Boris method \cite{boris70rps} applied to the non-relativistic equations of motion,
$$
\frac{x^{n+1}-2 x^n + x^{n-1}}{h^2} = E(x^n) + \frac{x^{n+1}-x^{n-1}}{2h}\times B(x^n) .
$$ 
For the energy-conserving leapfrog method the limit equation is similar, just with the discrete gradient $-\overline\nabla \phi^n$ instead of $E(x^n)=-\nabla \phi(x^n)$. This method preserves the energy $\tfrac12 |v^{n+1/2}|^2+ \phi(x^{n+1/2})$ for all $n$. 

For the variational leapfrog method the limit equation is
\begin{align*}
&\frac{x^{n+1}-2 x^n + x^{n-1}}{h^2} = E(x^n) + \frac{x^{n+1}-x^{n-1}}{2h}\times B(x^n) 
\\
&\qquad \qquad + A'(x^n) \frac{x^{n+1}-x^{n-1}}{2h} - \frac{A(x^{n+1})-A(x^{n-1})}{2h},
\end{align*}
which is the variational integrator for the non-relativistic equations of motion considered in \cite{hairer20lta,hairer22lsi,webb14sio}.

\section*{Acknowledgement} 
The work of Yanyan Shi was funded by the Sino-German (CSC-DAAD) Postdoc Scholarship, Program No. 57575640. Ernst Hairer acknowledges the support of the Swiss National Science Foundation, grant No.200020 192129.

\bibliographystyle{abbrv}
\bibliography{HLW}

\begin{thebibliography}{10}

\bibitem{benettin94oth}
G.~Benettin and A.~Giorgilli.
\newblock On the {H}amiltonian interpolation of near to the identity symplectic
  mappings with application to symplectic integration algorithms.
\newblock {\em J.\ Statist.\ Phys.}, 74:1117--1143, 1994.

\bibitem{boris70rps}
J.~P. Boris.
\newblock Relativistic plasma simulation-optimization of a hybrid code.
\newblock {\em Proceeding of Fourth Conference on Numerical Simulations of
  Plasmas}, pages 3--67, November 1970.

\bibitem{celledoni12per}
E.~Celledoni, V.~Grimm, R.~McLachlan, D.~McLaren, D.~O'Neale, B.~Owren, and
  G.~Quispel.
\newblock Preserving energy resp. dissipation in numerical {PDEs} using the
  {\textquotedblleft}average vector field{\textquotedblright} method.
\newblock {\em J. Comput. Phys.}, 231(20):6770--6789, 2012.

\bibitem{goldstein2011classical}
H.~Goldstein, C.~Poole, and J.~Safko.
\newblock {\em Classical Mechanics}.
\newblock Pearson Education India, 2011.

\bibitem{gonzalez96tia}
O.~Gonzalez.
\newblock Time integration and discrete {H}amiltonian systems.
\newblock {\em J. Nonlinear Sci.}, 6:449--467, 1996.

\bibitem{hairer18ebo}
E.~Hairer and C.~Lubich.
\newblock Energy behaviour of the {B}oris method for charged-particle dynamics.
\newblock {\em BIT}, 58:969--979, 2018.

\bibitem{hairer20lta}
E.~Hairer and C.~Lubich.
\newblock Long-term analysis of a variational integrator for charged-particle
  dynamics in a strong magnetic field.
\newblock {\em Numer. Math.}, 144(3):699--728, 2020.

\bibitem{hairer22lsi}
E.~Hairer, C.~Lubich, and Y.~Shi.
\newblock Large-stepsize integrators for charged-particle dynamics over
  multiple time scales.
\newblock {\em Numer. Math.}, 151(3):659--691, 2022.

\bibitem{hairer03gni}
E.~Hairer, C.~Lubich, and G.~Wanner.
\newblock Geometric numerical integration illustrated by the
  {S}t\"ormer--{V}erlet method.
\newblock {\em Acta Numerica}, 12:399--450, 2003.

\bibitem{hairer06gni}
E.~Hairer, C.~Lubich, and G.~Wanner.
\newblock {\em Geometric Numerical Integration. {S}tructure-Preserving
  Algorithms for Ordinary Differential Equations}.
\newblock Springer Series in Computational Mathematics 31. Springer-Verlag,
  Berlin, 2nd edition, 2006.

\bibitem{he16hov}
Y.~He, Y.~Sun, R.~Zhang, Y.~Wang, J.~Liu, and H.~Qin.
\newblock High order volume-preserving algorithms for relativistic charged
  particles in general electromagnetic fields.
\newblock {\em Physics of Plasmas}, 23(9):092109, 2016.

\bibitem{higuera2017structure}
A.~V. Higuera and J.~R. Cary.
\newblock Structure-preserving second-order integration of relativistic charged
  particle trajectories in electromagnetic fields.
\newblock {\em Physics of Plasmas}, 24(5):052104, 2017.

\bibitem{jackson1999classical}
J.~D. Jackson.
\newblock {\em Classical Electrodynamics}.
\newblock Wiley, Singapore, 3rd edition, 1999.

\bibitem{maeda82lfo}
S.~Maeda.
\newblock Lagrangian formulation of discrete systems and concept of difference
  space.
\newblock {\em Math. Japonica}, 27:345--356, 1982.

\bibitem{matsuyama17hoi}
A.~Matsuyama and M.~Furukawa.
\newblock High-order integration scheme for relativistic charged particle
  motion in magnetized plasmas with volume preserving properties.
\newblock {\em Computer Physics Comm.}, 220:285--296, 2017.

\bibitem{ripperda2018comprehensive}
B.~Ripperda, F.~Bacchini, J.~Teunissen, C.~Xia, O.~Porth, L.~Sironi,
  G.~Lapenta, and R.~Keppens.
\newblock A comprehensive comparison of relativistic particle integrators.
\newblock {\em The Astrophysical Journal Supplement Series}, 235(1):21, 2018.

\bibitem{vay2008simulation}
J.-L. Vay.
\newblock Simulation of beams or plasmas crossing at relativistic velocity.
\newblock {\em Physics of Plasmas}, 15(5):056701, 2008.

\bibitem{verlet67ceo}
L.~Verlet.
\newblock Computer ``experiments'' on classical fluids. {I}. {T}hermodynamical
  properties of {L}ennard-{J}ones molecules.
\newblock {\em Physical Review}, 159:98--103, 1967.

\bibitem{wang2021high}
Y.~Wang, J.~Liu, and Y.~He.
\newblock High order explicit {L}orentz invariant volume-preserving algorithms
  for relativistic dynamics of charged particles.
\newblock {\em J. Comput. Phys.}, 439:110383, 2021.

\bibitem{wang2016lorentz}
Y.~Wang, J.~Liu, and H.~Qin.
\newblock Lorentz covariant canonical symplectic algorithms for dynamics of
  charged particles.
\newblock {\em Physics of Plasmas}, 23(12):122513, 2016.

\bibitem{webb14sio}
S.~D. Webb.
\newblock Symplectic integration of magnetic systems.
\newblock {\em J. Comput. Phys.}, 270:570--576, 2014.

\bibitem{xiao2019explicit}
J.~Xiao and H.~Qin.
\newblock Explicit high-order gauge-independent symplectic algorithms for
  relativistic charged particle dynamics.
\newblock {\em Computer Physics Comm.}, 241:19--27, 2019.

\bibitem{zhang2015volume}
R.~Zhang, J.~Liu, H.~Qin, Y.~Wang, Y.~He, and Y.~Sun.
\newblock Volume-preserving algorithm for secular relativistic dynamics of
  charged particles.
\newblock {\em Physics of Plasmas}, 22(4):044501, 2015.

\end{thebibliography}
%\bibliography{/Users/hairer/bin/HLW}

\end{document}